\documentclass[a4paper,fleqn,12pt,twoside]{article}

\input{extern/layout_commands}
\usepackage{amsmath}
\usepackage{amssymb}

\newcounter{mynumcounter}
	       {%
		 \begin{list}{(\roman{mynumcounter})\hspace*{\fill}}%
		   {
		     \setlength{\topsep}{0cm}
		     \setlength{\partopsep}{0cm}
		     \setlength{\itemsep}{0ex}
		     \setlength{\parsep}{0cm}
		     \setlength{\leftmargin}{0cm}
		     \setlength{\itemindent}{8mm}		   
		     \setlength{\labelsep}{3mm}
		     \setlength{\labelwidth}{5mm}
		     \usecounter{mynumcounter}
		   }%
	       }
	       {\end{list}}

\newcommand{\eps}{\varepsilon}
\newcommand{\del}{\partial}

\newcommand{\dd}[2]{\frac{\del #1}{\del #2}}
\newcommand{\ddeval}[3]{\left.\dd{#1}{#2}\right|_{#3}}

\newcommand{\DD}[2]{\frac{d #1}{d #2}}
\newcommand{\DDeval}[3]{\left.\DD{#1}{#2}\right|_{#3}}

\newcommand{\Vol}{\operatorname{Vol}}

\newcommand{\graph}{\operatorname{graph}}

\newcommand{\tr}{\operatorname{tr}}
\renewcommand{\div}{\operatorname{div}}

\newcommand{\IR}{\mathbf{R}}

\newcommand{\CH}{\mathcal{H}}

\newcommand{\CM}{\mathcal{M}}

\newcommand{\id}{\operatorname{id}}

\newcommand{\rmd}{\mathrm{d}}

\newcommand{\const}{\mathrm{const}}

\newcommand{\dmu}{\,\rmd\mu}
\newcommand{\im}{\mathrm{im}}

\newcommand{\inj}{\mathrm{inj}}

\newcommand{\Connection}[1]{\smash{\sideset{^{#1}}{}{\mathop\nabla\nolimits}}}

\newcommand{\nabM}{\!\Connection{M}}

\newcommand{\Riemann}[1]{\smash{\sideset{^{#1}}{}{\mathop\mathrm{Rm}\nolimits}}}

\newcommand{\RiemM}{\!\Riemann{M}}

\newcommand{\Scalarcurv}[1]{\smash{\sideset{^{#1}}{}{\mathop\mathrm{Sc}\nolimits}}}

\newcommand{\ScalM}{\!\Scalarcurv{M}}
\newcommand{\ScalSig}{\!\Scalarcurv{\Sigma}}

\newcommand{\Acirc}{\hspace*{4pt}\raisebox{8.5pt}{\makebox[-4pt][l]{$\scriptstyle\circ$}}A}

\newcommand{\Acircbar}%
{\hspace*{4pt}\raisebox{8.5pt}{\makebox[-4pt][l]{$\scriptstyle\circ$}}\bar A}
\newcommand{\Kcircbar}%
{\hspace*{4pt}\raisebox{8.5pt}{\makebox[-4pt][l]{$\scriptstyle\circ$}}\bar K}

\newcommand{\half}{\tfrac{1}{2}}

\input{extern/draft_commands}
\newcommand{\Hmax}{{H_\mathrm{max}}}
\renewcommand{\span}{\mathrm{span}}
\newcommand{\diam}{\mathrm{diam}}

\newcommand{\IN}{\mathbf{N}}

\title{\titlefamily\Huge Surfaces with maximal constant mean curvature}
\author{
  {\titlefamily Jan Metzger}
  \\
  {\titlefamily\small jan.metzger@aei.mpg.de}
  \\
  {
    \titlefamily\small
    Albert-Einstein-Institut,
    Am M\"uhlenberg 1,
    D-14476 Potsdam,
    Germany.
  }
}
\date{}

\begin{document}
\hyphenation{}
\pagestyle{footnumber}
\maketitle
\thispagestyle{footnumber}
\begin{abst}%
  In this note we consider asymptotically flat manifolds with
  non-negative scalar curvature and an inner boundary which is an
  outermost minimal surface. We show that there exists an upper bound
  on the mean curvature of a constant mean curvature surface
  homologous to a subset of the interior boundary components. This
  bound allows us to find a maximizer for the constant mean curvature
  of a surface homologous to the inner boundary.

  With this maximizer at hand, we can construct an increasing family
  of sets with boundaries of increasing constant mean curvature. We
  interpret this familiy as a weak version of a CMC foliation.
\end{abst}
\section{Introduction}
\label{sec:introduction}
Consider a non compact three dimensional Riemannian manifold $(M,g)$
with compact interior boundary $\del M$, which is the only minimal
surface in $(M,g)$.  In this paper we investigate how large the mean
curvature $H$ of embedded, constant mean curvature (CMC) surface in
the homology class of $\del M$ can be. The main result in this paper
is an upper bound for this curvature. Combined with an area estimate,
we then show that there exists a CMC surface which attains this maximum.

The main motivation for this work are CMC foliations. These foliations
have been used successfully in general relativity to study the center
of mass of isolated systems \cite{Huisken-Yau:1996} and the Riemannian
Penrose inequality \cite{Bray:1997}. The existence result in
\cite{Huisken-Yau:1996} constructs a CMC foliation in the asymptotic
region near infinity. The natural question arises how far to the
interior these foliations can be extended. It is clear that in general
topological reasons imply non-existence of an entire smooth
foliation. This calls for a weak version of a CMC foliation.

Let us consider a different perspective. If the interior boundary
$\del M$ is an outermost minimal surface, that is $(M,g)$ does not
contain any other minimal surface, then it is straightforward to
construct a local CMC foliation near $\del M$,
cf.~lemma~\ref{thm:inner_foliation}. So another question is how far
this interior foliation can be extended outward, away from $\del M$. This
is by far an easier question than extending the foliation inward.

The reason is the following. Roughly speaking, if we consider a
potential CMC foliation reaching from $\del M$ to infinity, then the
mean curvature has to increase near $\del M$, and decrease as in
Euclidean space when approaching infinity, as the surfaces of the
foliation enclose increasing volume with $\del M$. This has two
implications. First, there is a maximal value of CMC along this
foliation, and second that there are two types of behavior. The first
type is portions along which CMC increases and the other is where CMC
decreases. The former includes the region near $\del M$ and the latter
the asymptotic region.

The region in which CMC increases is easier to handle, as the maximum
principle implies that the CMC surfaces along the foliation can not
touch. In the exterior region it is a lot harder to get control on the
separation of the surfaces, as one can see from the many different
foliations by spheres that are possible in $\IR^3$.

This result gives a partial answer to the above questions. In
section~\ref{sec:an-upper-bound} we show that there is a bound for the
maximal CMC of a surface homologous to $\del M$. This needs a lower
bound on the scalar curvature of $M$, $\ScalM\geq -C$, and uses the
fact that $\del M$ is area minimizing. The condition that the surface
be homologous to $\del M$ is necessary, as near maxima of the
sectional curvature a CMC foliation exits where homologically trivial
spheres have unbounded CMC \cite{ye:1991}.

In section~\ref{sec:existence} the curvature bound is used to
construct a surface which realizes the maximal CMC $\Hmax$. This
existence result needs the stronger assumption that $\ScalM\geq 0$, in
order to ensure that the area of CMC surfaces is bounded. A unique
such surface can be selected by demanding that it be the innermost
such surface.

We show that this surface bounds a region with $\del M$, which can be
regarded as a manifold with boundary, cf.~the discussion at the end of
section~\ref{sec:existence} and in section~\ref{sec:weak-cmc}. Using
the outer boundary as barrier, we can construct an increasing family
of sets bounded by surfaces with CMC ranging form $0$ at the horizon
to $\Hmax$ at the boundary.  This increasing family is a candidate for
a weak version of a CMC foliation reaching up to the
$\Hmax$-surface. We explore some basic properties in the second half
of section~\ref{sec:weak-cmc}.


%
\section{Preliminaries}
\label{sec:preliminaries}
We consider asymptotically flat manifolds $(M,g)$ with an inner
boundary $\del M$ which is an outermost minimal surface. Such
manifolds $M$ are called \emph{exterior regions}. The requirement of
asymptotic flatness means that there exists a compact set $K\subset M$
and a diffeomorphism $x:M\setminus K\to \IR^3\setminus B_1(0)$ such
that in the $x$-coordinates the metric $g$ approaches the Euclidean
metric $\delta$, that is there exists $C$ such that
\begin{equation*}
  r|g -\delta | + r^2|\del g| \leq C.
\end{equation*}
To say that $\del M$ is an outermost minimal surface means that there
does not exist another minimal surface in $M$ which is homologous to
$\del M$. For an asymptotically flat manifold which contains minimal
surfaces, the outermost minimal surface always exists and is unique
\cite[Section 4]{Huisken-Ilmanen:2001}. An exterior region $M$ is
diffeomorphic to $\IR^3\setminus (\bigcup_{i=1}^N B_i)$, where the
$B_i$ are open balls with disjoint closure. Hence $\del M =
\bigcup_{i=1}^N S_i$, where $S_i =\del B_i$. This restricted topology
does not require any curvature assumptions. The fact that $\del M$ is
an outermost minimal surface implies furthermore that for each
$I\subset \{1,\ldots,N\}$ the set $\bigcup_{i\in I} S_i$ minimizes
area in its homology class, in particular $\del M$ is minimizing.

Let $\Sigma\subset M$ be a two-sided surface. We assume that we can
identify one side of $\Sigma$ as the outside, and denote the outward
pointing normal by $\nu$. We denote by $\gamma$ the induced
metric. The mean curvature $H =\div \nu$ is taken with respect to the
outward pointing normal as is the second fundamental form $A$. By
$\ScalSig$ we denote the scalar curvature of $\Sigma$. The trace free
part of the second fundamental form will be denoted by $\Acirc = A -
\half H \gamma$.

Consider a normal variation of $\Sigma$, that is a map
$F:\Sigma\times(-\eps,\eps)\to M$ with $F(\cdot,0) =\id_\Sigma$ and
$\DDeval{F}{t} = f\nu$. The linearization of the operator which
assigns the mean curvature to the surfaces $\Sigma_t = F(\Sigma,t)$ is
given by
\begin{equation*}
  \ddeval{}{t}{t=0} F_t^* H(\Sigma_t)
  =
  Lf
  =
  -\Delta f - (\half \ScalM - \half\ScalSig + |\Acirc|^2 +
  \tfrac{3}{4} H^2)f
\end{equation*}
where $F_t = F(\cdot,t) : \Sigma\to\Sigma_t$ and $\Delta$ denotes the
Laplace-Beltrami operator along $\Sigma$. Here $\ScalM$ and $\ScalSig$
denote the scalar curvature of $M$ and $\Sigma$. $L$ is called the
\emph{stability operator}, or Jacobi operator.

When dealing with constant mean curvature surfaces $H=\const$, there
are two types of stability discussed in the literature. The first
notion is \emph{strong stability}, where we assume that $L$ is a
non-negative operator, that is
\begin{equation}
  \label{eq:1}
  \int_\Sigma
  f^2 (\half \ScalM - \half\ScalSig + \half|\Acirc|^2 + \tfrac{3}{4} H^2)
  \dmu
  \leq
  \int_\Sigma |\nabla f|^2 \dmu
  \qquad
  \forall f\in C^\infty(\Sigma).
\end{equation}
Here $\nabla f$ denotes the tangential gradient of $f$. Note that
strong stability means that the principal eigenvalue of $L$, that is
the smallest eigenvalue, is non-negative. The second notion, simply
called \emph{stability} comes from the fact that the constant mean
curvature equation is the Euler-Lagrange equation for the
isoperimetric problem, that is for minimizing the area of $\Sigma$,
while keeping enclosed volume constant. Minimizers of this variational
principle satisfy the stability inequality
\begin{multline*}
  \int_\Sigma
  f^2 (\half \ScalM - \half\ScalSig + \half|\Acirc|^2 + \tfrac{3}{4} H^2)
  \dmu
  \leq
  \int_\Sigma |\nabla f|^2 \dmu
  \\
  \forall f\in C^\infty(\Sigma)
  \ \text{with}\
  \int_\Sigma f \dmu
  =
  0.
\end{multline*}
Hence, strong stability implies stability, but not vice-versa. For the
following discussion only strong stability plays a role.

Subsequently, we deal with surfaces which are not necessarily
connected. We say that such a surface is strongly stable if each of
its components is strongly stable, and thus if a surface is not
strongly stable, it means that at least one of its components is not
strongly stable.

The surfaces $\Sigma$ in question will be homologous to $\del M$. In
the case that $\Sigma$ does not touch $\del M$ this means that there
exists an open set $\Omega$ such that $\del\Omega$ is the disjoint
union $\del\Omega = \del M \cup \Sigma$. As we orient $\del M$ with
the normal pointing into $M$, the correct orientation of $\Sigma$
corresponds to the normal vector pointing out of $\Omega$. We will
make this assumption subsequently without further notice.


%
\section{An upper bound for CMC}
\label{sec:an-upper-bound}
This section is devoted to derive an upper bound for the constant mean
curvature of a compact, smooth, embedded CMC surface homologous to
$\del M$. Note that $\Sigma$ need not be connected for the subsequent
arguments. This upper bound only requires a lower bound on the scalar
curvature of $M$, that is $\ScalM\geq -C$ for some $C\geq 0$.

Before we can approach the main theorem, we review an existence
theorem for prescribed mean curvature surfaces \cite[Theorem
6.1]{Andersson-Metzger:2007}. This theorem implies the following
existence theorem for strongly stable CMC surfaces.
\begin{theorem}
\label{thm:exist_barrier}
Let $(\Omega,g)$ be a compact Riemannian manifold with smooth boundary
$\del\Omega$ which is the disjoint union
$\del\Omega=\del^-\Omega\cup\del^+\Omega$, where $\del^\pm\Omega$ are
smooth, non-empty and without boundary. Assume that $\del^-\Omega$ has
mean curvature $H^-$, where $H^-$ is taken with respect to the normal
pointing into $\Omega$, and $\del^+\Omega$ has mean curvature $H^+$,
where $H^+$ is taken with respect to the normal pointing out of
$\Omega$. Let $h$ be such that $\max_{\del^-\Omega} H^- \leq h \leq
\min_{\del^+\Omega} H^+$. Then there exists a compact, smooth,
embedded, strongly stable CMC surface $\Sigma\subset \Omega$,
homologous to $\del^-\Omega$ with $H(\Sigma) = h$.
\end{theorem}
\begin{proof}
  The theorem is a direct consequence of Theorem 6.1 from
  \cite{Andersson-Metzger:2007}. This theorem states that if $(\Omega,
  g)$ is as in the assumption and $K$ is a symmetric bilinear form on
  $\Omega$, then if $\theta^+(\del^-\Omega) \leq 0$ and
  $\theta^+(\del^+\Omega)\geq 0$ then there exists a compact, smooth,
  embedded, surface $\Sigma$ homologous to $\del^-\Omega$ with
  $\theta^+(\Sigma)=0$, which is stable in the sense of surfaces with
  $\theta^+=0$. Here $\theta^+ = H + P$, where $H$ is the mean
  curvature as usual, and $P = \tr_\Sigma K = \tr_M K - K(\nu,\nu)$ is
  the trace of $K$ restricted to $T\Sigma$.

  We apply this theorem to the data $(\Omega,g, K = - \half h g)$ such
  that for any surface $\Sigma$ we have $\theta^+(\Sigma) = H -
  h$. Thus $\theta^+(\del^-\Omega) = H^- -h \leq 0$ and
  $\theta^+(\del^+\Omega) = H^+ - h \geq 0$, and the existence of a
  surface $\Sigma$ with $H(\Sigma) = h$ follows from the existence
  part of theorem 6.1 in \cite{Andersson-Metzger:2007}.  The resulting
  surface is also stable in the sense of $\theta^+=0$ surfaces, which
  implies that the smallest eigenvalue of the operator
  \begin{equation*}
    \tilde L f
    =
    - \Delta f
    + 2 S(\nabla f)
    - f (\div S  - \half |A + K^\Sigma|^2  - |S|^2  + \half \ScalSig -
    \mu + J(\nu) ) 
  \end{equation*}
  is non-negative. Here $K^\Sigma = K|_\Sigma$, $S = K(\nu,\cdot)^T$,
  where $T$ denotes tangential projection, $\mu = \half (\ScalM -
  |K|^2 + (\tr_M K)^2)$ and $J = \div K -\nabM \tr_M K$. On a surface
  with $\theta^+=0$ we have for our choice of $K$ that $S =0$, $J=0$,
  $|A + K^\Sigma|^2 = |A|^2 - \half H^2 = |\Acirc|^2$ and $\mu = \half
  \ScalM + \tfrac{3}{4} H^2$. Thus we find that $\tilde L$ is nothing
  but the stability operator $L$ and non-negativity of its first
  eigenvalue means strong stability.
\end{proof}
\begin{remark}
  A similar existence theorem could be derived by analyzing the
  functional 
  \begin{equation*}
    J_h(F) = |\del^* F| - h \Vol(F)
  \end{equation*}
  for sets $F$ with finite perimeter.
\end{remark}  
Using the previous existence theorem together with the fact that $\del
M$ can be used as inner barrier, we infer the following lemma.
\begin{lemma}
  \label{thm:exist-strongly-stable}
  If $\Sigma\subset M$ is a CMC surface with $H>0$ homologous to $\del
  M$, then there exists a strongly stable CMC surface $\Sigma'$ in the
  same homology class with $H(\Sigma') = H(\Sigma)$.
\end{lemma}
\begin{proof}
  \begin{figure}[!h]
    \centering
    \resizebox{.5\linewidth}{!}{\input{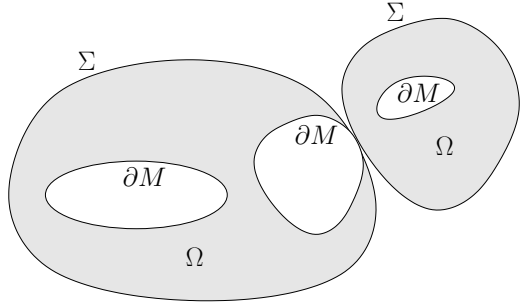}}
    \caption{If a surface $\Sigma$ homologous to $\del M$ intersects
      $\del M$, then there is one component of $\Sigma$ which intersects
      $\del M$ such that the outer normals point in the same
      direction.}
    \label{fig:touch}
  \end{figure}
  First, note that $\Sigma$ can not touch $\del M$, cf.\ figure
  \ref{fig:touch}. As $\Sigma$ is homologous to $\del M$, there exists
  a set $\Omega$ and a set $I \subset \{1,\ldots,N\}$ such that $\del
  \Omega = \bigcup_{i\in I} S_i \cup \Sigma_0$, and $\Sigma = \Sigma_0
  \cup \bigcup_{i\not\in I} S_I$. As $H(\Sigma)>0$, we have must have
  $I=\{1,\ldots,N\}$, and thus $\Sigma_0 = \Sigma$ and $\del M\subset
  \Omega$. Thus if $\del M$ and $\Sigma$ intersect, there exists a
  component $\Sigma_1$ of $\Sigma$ which intersects $\del M$ at a
  point where the normals of $\Sigma_1$ and $\del M$ point in the same
  direction. This is impossible, since the maximum principle would
  imply that $\Sigma_1\subset \del M$.

  Hence, $\Sigma$ lies completely in the interior of $M$ and we can
  apply theorem~\ref{thm:exist_barrier} with the so constructed
  $\Omega$ where $\del M = \del^-\Omega$, $\Sigma = \del^+\Omega$ and
  $h = H(\Sigma)$. Thus we obtain $\Sigma'$, a smooth, embedded
  strongly stable CMC surface.  
\end{proof}
\begin{remark}
  Note that we can in fact show that the constructed $\Sigma'$ can not
  touch any component of $\Sigma$ which is not strongly stable, as
  these components can be deformed in direction of $-\nu$, that is
  into $\Omega$ in such a way that their mean curvature increases.
\end{remark}
It is now a simple matter to derive the claimed bound on the constant
mean curvature from strong stability. As $\Sigma$ is not necessarily
connected, we must make use of the fact that $\del M$ is outermost to
get a lower bound on the area of at least one component of $\Sigma$.
\begin{lemma}
  \label{thm:supbound}
  Let $(M,g)$ be asymptotically flat with inner boundary $\del M$,
  which is an outermost minimal surface in $M$. Assume that
  \begin{equation*}
    \ScalM\geq -C
  \end{equation*}
  Denote the components of $\del M$ by $S_i$,
  $i=1,\ldots,N$ and let $A:=\max \{ |S_i| : i=1,\ldots,N\}$. If
  $\Sigma\subset M$ is a CMC surface homologous to $\del M$, then
  \begin{equation*}
    H(\Sigma)^2 \leq \frac{16\pi}{3A}  + \frac{2}{3}C.
  \end{equation*}
\end{lemma}
\begin{remark}
  An obvious modification yields a similar bound if $\Sigma$ is
  homologous to $\bigcup_{i\in I} S_i$, where $I\subset
  \{1,\ldots,N\}$ and $A$ is replaced by $A(I) = \max_{i\in I} |S_i|$.
\end{remark}
\begin{proof}
  Assume that $A = |S_1|$. As $M$ is topologically equivalent to
  $\IR^3 \setminus \bigcup_{i=1}^N B_i$, as explained in
  section~\ref{sec:preliminaries}, $\Sigma$ can be regarded as a
  surface embedded in $\IR^3$. Then any component of $\Sigma$ bounds in
  $\IR^3$, and since $\Sigma$ is homologous to $\bigcup_{i=1}^N S_i$,
  we infer that there exists one component $\Sigma_1$ of $\Sigma$
  which is homologous to $S_1 \cup \bigcup_{i\in J} S_i$, where $J
  \subset \{2,\ldots,N\}$ may be empty. Since $S_1 \cup \bigcup_{i\in
    J}S_i$ is minimizing in its homology class in $M$ we find that
  $|\Sigma_1| \geq |S_1 \cup \bigcup_{i\in J}S_i| \geq |\Sigma_1| =
  A$. Pick a test function $f\in C^\infty(\Sigma)$ with $f=1$ on
  $\Sigma_1$ and $f=0$ on all other components. Plugging $f$ into the
  strong stability inequality~\eqref{eq:1}, we find that
  \begin{equation*}
    \int_{\Sigma_1}
    \half|\Acirc|^2 + \tfrac{3}{4} H^2
    \dmu
    \leq
    \int_{\Sigma_1} \half\ScalSig -\half\ScalM \dmu.
  \end{equation*}
  From Gauss-Bonnet we infer that
  \begin{equation*}
    \int_{\Sigma_1} \half\ScalSig\dmu = 4\pi(1-\mathrm{genus}(\Sigma_1)) \leq 4\pi.
  \end{equation*}
  As $H$ is constant, combining the above inequality with the lower
  bound on $\ScalM$ yields
  \begin{equation}
    \label{eq:3}
    H^2|\Sigma_1| \leq \frac{16\pi}{3} + \frac{2}{3} C|\Sigma_1|,
  \end{equation}
  or after dividing by $|\Sigma_1|$,
  \begin{equation*}
    H^2 \leq \frac{16\pi}{3|\Sigma_1|} + \frac{2}{3} C,
  \end{equation*}
  which implies the claim, since $|\Sigma_1|\geq A$.
\end{proof}
\begin{remark}
  The spatial Schwarzschild manifold of mass $m$ is $(\IR^3\setminus
  \{0\}, \phi^4 g^e)$ where $\phi= 1 +\frac{m}{2r}$ and $g^e$ denotes
  the Euclidean metric on $\IR^3$. It is scalar flat and if $m>0$ it
  has an outermost minimal surface at $r = \frac{m}{2}$. Thus
  $(\IR^3\setminus B_{\frac{m}{2}}, \phi^4 g^e)$ satisfies the
  assumptions of lemma~\ref{thm:supbound} with $C=0$. The spheres
  $S_r(0)$ have constant mean curvature $H_r = \frac{2}{R}
  \frac{2r-m}{2r+m}$ where $R=\phi^2 r$ is the geometric area radius
  of $S_r$ with respect to $g^S$, that is $|S_r| = 4\pi R^2$. $H_r$
  assumes its maximum where $R=3m$ and equals $\frac{2}{3\sqrt{3}m}$
  there. Thus, the estimate of equation~\eqref{eq:3} is sharp in this
  case, whereas the assertion of lemma~\ref{thm:supbound} is not.
\end{remark}


%
\section{Existence of  surfaces with maximal CMC}
\label{sec:existence}
In this section we construct a surface with maximal constant mean
curvature. In fact, for $(M,g)$ as before, we can let
\begin{equation*}
  \Hmax := \sup \{ H(\Sigma) : \Sigma\ \text{an embedded CMC surface
    homologous to}\ \del M\}.
\end{equation*}
As we have seen in the previous section, $\Hmax$ is finite.
Subsequently we show that $\Hmax$ is attained at a strongly
stable surface. We start by showing that $\Hmax>0$.
\begin{lemma}
  \label{thm:inner_foliation}
  There exists a foliation of a neighborhood of $\del M$ by CMC
  surfaces $\Gamma_s$, $s\in[0,\eps)$ with $H(\Gamma_s) >0$.
\end{lemma}
\begin{proof}
  We construct the foliation near each component of $\del M$
  separately. Let $S_i$ be such a component. Note that $S_i$ is stable
  as a minimal surface, as $\del M$ is outermost. If the principal
  eigenvalue $\lambda$ of $L$ on $S_i$ is positive, $\lambda>0$, then
  $L$ is invertible and we can construct a foliation of CMC surfaces
  by a simple application of the implicit function theorem.

  Hence we can assume $\lambda=0$ from now on. Let $\phi>0$ denote the
  corresponding eigenfunction. In this case a CMC foliation can be
  constructed as in \cite{Galloway:2006ws}. We repeat the argument
  here for convenience. Consider the operator
  \begin{equation*}
    \CH
    :
    C^\infty(\Sigma) \times \IR \to C^\infty \times \IR
    :
    (u,h) \mapsto \left( H(\graph u) - h, \int_\Sigma u\phi \dmu\right),
  \end{equation*}
  where $\graph u = F_u(\Sigma)$ and $F_u(p) = \exp_p(u(p)\nu_p)$,
  where $p\in\Sigma$ and $\exp$ is the exponential map of $M$. Then
  $H(\graph u)$ denotes the mean curvature of $\graph u$ pulled-back
  to $\Sigma$ via $F_u$.

  We can compute the linearization of $\CH$ at $(u,h) = (0,0)$ in
  direction $(v,s)\in C^\infty(\Sigma)\times\IR$ to be
  \begin{equation*}
    \CM := D\CH |_{(0,0)} (v,s) = \left(Lv - s, \int_\Sigma v\phi \dmu\right).
  \end{equation*}
  Obviously $\CM$ is invertible since $\ker L = \span \{\phi\}$ and
  the equation $Lv = g$ is uniquely solvable if $\int_\Sigma g\phi
  \dmu = 0$ and $\int_\Sigma v \phi \dmu = 0$.

  By the inverse function theorem, there exists $u(t)$ and $h(t)$ for
  small $t$ such that
  \begin{equation}
    \label{eq:4}
    \CH(u(t),h(t)) = (0,t).
  \end{equation}
  This implies that the surfaces $\graph u(t)$ have
  CMC. Differentiating equation~\eqref{eq:4} with respect to $t$
  yields that
  \begin{equation}
    \label{eq:2}
    \left( L u'(0) - h'(0), \int_\Sigma u'(0) \phi \dmu \right)
    =
    \big(0,1\big)
  \end{equation}
  and hence that $h'(0) \in \im L \perp \ker L$, that is $\int h'(0)
  \phi \dmu =0$. Since $h'(0)$ is a constant and $\phi>0$, we infer
  $h'(0) =0$. Then $u'(0) \in\ker L$ and $u'(0) = \alpha\phi$ where
  $\alpha >0$, by~\eqref{eq:2}. Thus, the $\graph u(t)$ form a
  foliation near $S_i$.

  As $\del M$ is outermost, we must have that $h(t)>0$ for all $t$,
  and we thus found the foliation near $S_i$. As $h(t)$ is smooth,
  there exists a $t_0$ such that $h$ is increasing on $[0,t_0)$.

  Thus we can find the required CMC foliation near each component of
  $\del M$ separately and join it to give a CMC foliation near $\del
  M$.
\end{proof}
\begin{remark}
  A different way to see that $\Hmax>0$ is to use asymptotic flatness
  to conclude that there exists a surface in the asymptotic end with
  positive mean curvature. An application of
  theorem~\ref{thm:exist_barrier} then yields a CMC surface with
  positive CMC. However, the previous lemma emphasizes that not only
  the asymptotic behavior near infinity, but also the local geometry
  near $\del M$ gives a lower bound on $\Hmax$.
\end{remark}
Standard arguments show that there are uniform bounds on the second
fundamental form of strongly stable CMC surfaces.
\begin{lemma}
  \label{thm:curv_bound}
  If $\Sigma$ is a strongly stable CMC surface then there exists a
  constant $C= C(\|\RiemM\|_{C^0}, \inj(M,g)^{-1}, \sup_\Sigma |H|)$ such that
  \begin{equation*}
    \sup_\Sigma |A| \leq C.
  \end{equation*}
\end{lemma}
\begin{proof}
  First, there exists $0 < r_0 = r_0(\|\RiemM\|_{C^0},\sup_\Sigma
  |H|)$ such that for all $r<r_0$ and $p\in\Sigma$ the area of the
  intrinsic balls $B^\Sigma(p,r)$ around $p$ with radius $r$ is
  bounded
  \begin{equation*}
    |B^\Sigma(p,r)| \leq 6\pi r^2.
  \end{equation*}
  See for example \cite[Theorem 8.1]{Andersson-Metzger:2005}, which
  goes back to \cite{Pogorelov:1981}. With this local bound on area,
  the usual argument for deriving curvature bounds yields the desired
  estimate, cf.~\cite{Schoen-Simon-Yau:1975}, we refer to \cite[Section
  6]{Andersson-Metzger:2005} for a detailed proof in a slightly more
  general setting.
\end{proof}
Before we can attempt the construction of surfaces realizing $\Hmax$,
we need a diameter bound for strongly stable CMC surfaces
\cite{Rosenberg:2006}.
\begin{lemma}
  \label{thm:diameter_bound}
  Let $(M,g)$ be a complete Riemannian 3-manifold with $\ScalM\geq 0$
  and let $\Sigma \subset M$ be a closed, connected, strongly stable
  CMC surface with $H(\Sigma) \neq 0$. Then
  \begin{equation*}
    \diam (\Sigma) \leq \frac{2\pi}{3H}.
  \end{equation*}
\end{lemma}
\begin{proof}
  This estimate is a direct consequence of \cite[Theorem
  1]{Rosenberg:2006}.
\end{proof}
\begin{theorem}
\label{thm:existence}
Let $(M,g)$ be an asymptotically flat Riemannian manifold with
$\ScalM\geq 0$ and a non-empty inner boundary $\del M$, which is an
outermost minimal surface. Assume that $\|\RiemM\|_{C^0}$ is finite
and $\inj(M,g)$ is non-zero. Then $\Hmax$ is attained at a compact,
immersed, strongly stable surface $\Sigma$ homologous to $\del M$.
$\Sigma$ is a union of spheres.
\end{theorem}
\begin{proof}
  Let $\{\Sigma^n\}_{\{n\geq 1\}}$ be family of CMC surfaces
  homologous to $\del M$ with
  \begin{equation*}
    H(\Sigma^n)\to\Hmax.
  \end{equation*}
  We show that after suitable modification, the sequence $\Sigma^n$
  allows the extraction of a convergent subsequence.

  In view of lemma~\ref{thm:exist-strongly-stable} we can assume that
  the $\Sigma^n$ are strongly stable. Due to
  lemma~\ref{thm:inner_foliation} we can furthermore assume that
  $H(\Sigma^n)\geq \eps$ for some suitably chosen $\eps>0$. 

  Fix an arbitrary $n$ and denote $\Sigma := \Sigma^n$. Let $\Sigma_j$
  be the components of $\Sigma$, $j=1,\ldots N_\Sigma$. For
  $j=1,\ldots,N_\Sigma$ let $f_j$ be the test function which is equal
  to $1$ on the $\Sigma_j$ and $0$ on the other components, and plug
  $f_j$ in the strong stability inequality~\eqref{eq:1}. This yields
  that $\Sigma_j$ is a sphere, as $\int_{\Sigma_j} H(\Sigma)^2 \dmu >0$.

  Let
  \begin{equation*}
    J := \{ j : \Sigma_j\ \text{does not bound a compact region on its
      inside} \}
  \end{equation*}
  and delete all components $\Sigma_j$ from $\Sigma$ where $j
  \not\in J$. The surface
  \begin{equation*}
    \Sigma' := \bigcup_{j\in J} \Sigma_j
  \end{equation*}
  is homologous to $\del M$ and thus separates $\del M$ from infinity.
  Recall that $M$ is diffeomorphic to $\IR^3\setminus \bigcup_{i=1}^N
  B_i$ and consider $\Sigma'\subset \IR^3$.

  Let $U\subset\IR^3$ be such that $\IR^3\setminus U$ is the
  non-compact component of $\IR^3 \setminus \Sigma'$.  Note that $\del
  U$ consists of a subset of the components of $\Sigma'$. Indeed $\del
  U = \Sigma'$. Otherwise there exists one component $\Sigma_j$ not in
  $U$ which bounds a compact region $\Omega_j$ on its outside,
  relative to a subset of $\del M$. This is clearly impossible as the
  boudary of $M':=M\setminus \Omega_j$ has $H(\del M') \leq 0$ and
  $H(\del M') \neq 0$. This would imply the existence of a minimal
  surface outside of $\del M$, contradicting the assumption that $\del
  M$ is an outermost minimal surface.

  Then $\Sigma' := \del U$ has at most $N$ components, each of which
  is homologous to $\bigcup_{I_j} S_i$, where $I_j\subset
  \{1,\ldots,N\}$ is non-empty, and $I_j \cap I_{j'} = \emptyset$ for
  $j\neq j'$. To see this, let $U_j$ be the compact region in $\IR^3$
  bounded by $\Sigma_j$. Then $U_j$ contains at least one of the
  $B_i$, so it is clear that $I_j \neq \emptyset$. Since we have
  $\Sigma' = \del U$, no component of $\Sigma'$ is separated from
  infinity by another component of $\Sigma'$, in particular the outer
  normal direction of $\Sigma'$ agrees with the outer normal to $\del
  U$. As all the $B_i$ are contained in $U$ and the components of
  $\Sigma'$ can not intersect, this implies that the each $B_i$ can be
  in at most one $U_j$. Thus the $I_j$ are mutually disjoint. For
  subsequent use we relabel the $(\Sigma^n)'$ as $\Sigma^n$.
  
  We thus have constructed a sequence $\{\Sigma^n\}$ of CMC surfaces
  with $H(\Sigma^n) \to \Hmax$, where each of the $\Sigma^n$ is a
  strongly stable CMC surface with at most $N$ components, and each
  component is homologous to a non-empty union of components of $\del
  M$.

  As $H(\Sigma^n)\geq \eps$, lemma~\ref{thm:diameter_bound} implies
  that each component of $\Sigma^n$ has bounded diameter. Such a
  component of $\Sigma^n$ encloses at least one of the $S_i$. We thus
  infer that there exists a compact set $B \subset M$ such that
  $\Sigma^n \subset B$ for all $n$. Furthermore, the curvature
  estimates from lemma~\ref{thm:curv_bound} imply uniform curvature
  bounds for $\Sigma^n$. Therefore the Ricci curvature of $\Sigma_n$
  is bounded below and standard volume comparison shows that each
  component of $\Sigma^n$ has bounded area. As there are at most $N$
  components the $\Sigma^n$ have uniformly bounded area.
  
  These three estimates, area, curvature, and the fact that the
  $\Sigma_n$ are contained in a compact set, imply that there exists a
  convergent subsequence and a limiting surface $\Sigma$, which has
  CMC and consists of strongly stable components. Note that the limit
  $\Sigma$ might not be embedded. Nevertheless, $\Sigma$ has an
  outward pointing normal vector field $\nu$ which is the limit of
  the outward pointing normal vector fields of the subsequence of
  $\Sigma_n$.
\end{proof}
We now examine the limiting $\Hmax$-surface $\Sigma$ more
closely. As $\Sigma$ is the limit of embedded surfaces, $\Sigma$ can
fail to be embedded only if $\Sigma$ touches itself. Let $p\in M$
denote such a point. Then at $p$ multiple sheets of $\Sigma$ can come
together. Since we have bounded curvature and bounded area, there are
at most finitely many such sheets $\Sigma^p_k$, $k=1,\ldots n(p)$, as
each sheet takes up some area.

Around $p$ there are coordinates $\{x^i\}$ of $M$ such that the
$\Sigma^p_k$ are $C^\infty$ graphs over an open subset $U^p$ of the
$x^1,x^2$-plane. That is $\Sigma^p_k = \{ x: x^3 = u(x^1,x^2) \}$. We
can assume that $u_k \leq u_l$ whenever $k\leq l$ since $\Sigma$ is
the limit of embedded surfaces which bound a region with respect to
$\del M$. Each of the sheets comes with a normal vector field $\nu^p$
with respect to which $H =\Hmax$. This can be either the downward or
upward pointing normal to $\graph u_\tau$, and this direction
alternates.

\begin{figure}[!t]
  \begin{minipage}{.45\linewidth}
    \resizebox{\linewidth}{!}{\input{pics/touch_out_pic.tex}}
    \caption{Two sheets touching on the outside.}
    \label{fig:touch_out}
  \end{minipage}
  \hfill
  \begin{minipage}{.45\linewidth}
    \resizebox{\linewidth}{!}{\input{pics/touch_in_pic.tex}}
    \caption{Two sheets touching on the inside.}
    \label{fig:touch_in}
  \end{minipage} 
\end{figure}
We say that two sheets $\Sigma_1$ and $\Sigma_2$ touch on the
\emph{outside} at a point $p$, if the representing functions $u_1
\leq u_2$ of these sheets are so that the outward normal of $\Sigma$
points upward along $u_1$ and downward along $u_2$, cf.
figure~\ref{fig:touch_out}.
On the other hand, if the normal along $u_1$ points downward, and
upward along $u_2$, we say that $\Sigma$ touches itself on the
\emph{inside}.

The following theorem is a direct consequence of the strong maximum
principle for surfaces with prescribed mean curvature.
\begin{theorem}
  \label{thm:touch_outside}
  Let $\Sigma$ be the $\Hmax$ surface constructed in
  theorem~\ref{thm:existence}. Then if $\Sigma$ is not embedded,
  $\Sigma$ can only touch itself on the outside, and no more than two
  sheets of $\Sigma$ can meet at one point of $M$.
\end{theorem}
\begin{proof}
  Let $\Sigma_1$ and $\Sigma_2$ be two sheets of $\Sigma$ which meet on
  the inside, and let $u_1\leq u_2$ be the representing functions as
  described above. Instead of the upward normal, consider $u_2$
  equipped with the downward normal. Then the mean curvature of
  $\graph u_2$ is $-\Hmax<0$ with respect to the downward normal, and
  the mean curvature of $\graph u_1$ is $\Hmax$ with respect to the
  downward unit normal. As $\graph u_1$ and $\graph u_2$ touch, we
  immediately obtain a contradiction to the strong maximum principle.

  At any point in $M$ where three sheets of $\Sigma$ meet, $\Sigma$
  must touch itself on the inside, thus this is ruled out by the
  above argument.
\end{proof}
We now want to add a few remarks about the uniqueness of the $\Hmax$
surfaces. Indeed, we can single out one particular $\Hmax$-surface in
$(M,g)$ by choosing the \emph{innermost} $\Hmax$ surface. 
\begin{theorem}
  \label{thm:innermost}
  Let $(M,g)$ be as in theorem~\ref{thm:existence}. Then there exists
  a unique innermost surface in $(M,g)$ which is homologous to $\del
  M$ and has CMC $\Hmax$. The assertion of
  theorem~\ref{thm:touch_outside} holds for $\Sigma$.
\end{theorem}
\begin{proof}
  The construction of this surface is similar to the construction of
  the outermost MOTS in \cite[Section 7]{Andersson-Metzger:2007}. Thus
  we mention only the key points for the construction.

  \emph{Compactness.} As in the proof of theorem~\ref{thm:existence},
  we infer compactness of the class of $\Hmax$-surfaces by throwing
  away components which bound compact regions.

  \emph{Monotonicity.} Let $\Sigma_i$, with $i=1,2$ be two
  $\Hmax$-surfaces for which the assertion of
  theorem~\ref{thm:touch_outside} holds, and which are homologous to
  $\del M$ and bound sets $\Omega_i$ with $\del M$. Then $\Omega_1
  \cap \Omega_2$ contains a strongly stable $\Hmax$-surface $\Sigma$
  homologous to $\del M$ which satisfies the assertion of
  theorem~\ref{thm:touch_outside}.

  Monotonicity allows us to construct a sequence of surfaces $\Sigma_k$
  bounding $\Omega_k$ together with $\del M$, such that the $\Omega_k$
  are decending. By compactness we find a limiting set $\Omega_\infty$
  bounded by an $\Hmax$-surface $\Sigma_\infty$ and $\del M$. 
\end{proof}
  

%
\section{A proposal for a weak CMC foliation}
\label{sec:weak-cmc}
In this section we propose a weak version of a foliation by CMC
surfaces of the interior region of $(M,g)$. There is more than one way
to introduce such a foliation, and it is not clear whether the
possibility discussed below is best suited for applications.

Let $(M,g)$ be asymptotically flat with $\del M$ an outermost minimal
surface. Assume that $\ScalM\geq 0$ and let $\Sigma$ be the
$\Hmax$-surface homologous to $\del M$ constructed in
section~\ref{sec:existence}.
\begin{figure}[!h]
  \centering
  \resizebox{.45\linewidth}{!}{\input{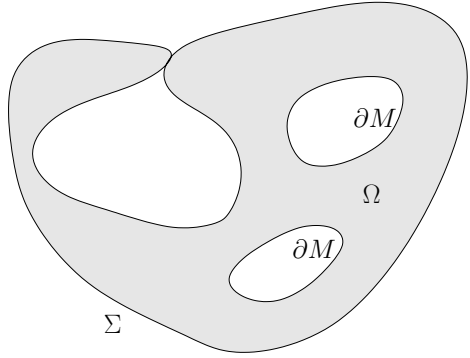}}
  \caption{The interior region $\Omega$.}
  \label{fig:interior}
\end{figure}
We denote by the interior $\Omega$ of $M$ the components of
$M\setminus\Sigma$ which meet components of $\del M$, cf.\ figure
\ref{fig:interior}. As $\Sigma$ does not touch itself on the inside,
$\Omega$ can be equipped with the structure of a smooth manifold with
boundary $\del M \cup \Sigma$, where we identify $\Sigma$ and $\del M$
with the points added by the metric completion of $\Omega$. In this
way, we separate the points of $\Sigma$ which are mapped to the same
point in the immersion of $\Sigma$ into $M$. Note that the interior
$\Omega$ is not a submanifold with smooth boundary in $M$ if $\Sigma$
is not embedded. The boundary of $\Omega$ consists of $\del M$ on the
inside, subsequently denoted by $\del^-\Omega$, and $\Sigma$ on the
outside, subsequently denoted by $\del^+\Omega$.

\subsection{Construction}
To construct a weak CMC foliation for this new manifold $\Omega$,
we introduce the following notion.
\begin{definition}
  \label{def:outermost}
  Let $\Sigma\subset \Omega$ be a smooth, embedded surface homologous
  to $\del^-\Omega$, with CMC $h\in(0,\Hmax)$. Denote by $U$ the
  region bounded by $\Sigma$ and $\del^+\Omega$. If there does not
  exist a smooth embedded surface $\Sigma'$ in $U$ with the same CMC
  $h$, then $\Sigma_h$ is called \emph{outermost}.
\end{definition}
By \cite[Section 7]{Andersson-Metzger:2007}, for each $h\in(0,\Hmax)$
there exists a smooth, embedded surface $\Sigma_h$, homologous to
$\del^-\Omega$, which has CMC $h$ and is outermost in the sense of
definition~\ref{def:outermost}. We denote by $\Omega_h$ the open
region bounded by $\del^-\Omega$ and $\Sigma$. We define this family
of sets to be the candidate for our weak CMC foliation of $\Omega$.

A useful side-effect of this definition is that the constructed sets
are related to a variational principle. Consider sets $F$ of finite
perimeter in $\Omega$. We will assume that $F\supset\Omega_h$ for some
$h>0$. Thus $F$ has one boundary component which agrees with
$\del^-\Omega$. In accordance with the above notation, we denote by
$\del^+ F$ the reduced boundary of $F$ without $\del^-\Omega$, that is
$\del^+ F = \del^*F \cap \Omega$.

\subsection{Basic properties}
Consider the functional $J_h$, defined on the collection of sets $F$
of bounded perimeter in $\Omega$,
\begin{equation*}
  J_h(F) := |\del^+ F| - h \Vol(F).
\end{equation*}
The critical points of $J_h$ are surfaces with CMC $h$, so it is
natural to consider this functional here.

We say that a set $E$ minimizes $J_h$ on the outside, if for all sets
$F\supset E$ we have
\begin{equation*}
  J_h(E) \leq J_h(F).
\end{equation*}
\begin{lemma}
  \label{thm:outward_Jh}
  For each $h\in(0,\Hmax)$, the set $\Omega_h$ defined above minimizes
  $J_h$ on the outside.
\end{lemma}
\begin{proof}
  If $\Omega_h$ does not minimize $J_h$ on the outside, then there
  exits a minimizer $E_h$ for $J_h$ outside of $\Omega_h$, with $E_h
  \neq \Omega_h$. The outer boundary $E_h$ is a
  $C^{1,\alpha}$-surface, satisfying $H\geq h$ in a distributional
  sense, cf.~\cite[Theorem 1.3]{Huisken-Ilmanen:2001}. It is smooth
  with $H=h$ where it does not touch $\Sigma_h$. By the strong maximum
  principle (which applies here as $\del^+ E_h$ is $C^{1,\alpha}$, see
  also \cite[Section 3]{Ziemer-Zumbrum:1998}) all components of
  $\del^+ E_h$ touching $\Sigma_h$ are contained in $\Sigma_h$. Thus,
  $\del^+ E_h$ is a smooth surface with CMC $h$ and lies on the
  outside of $\Sigma_h$. As $\Sigma_h$ is outermost, $E_h =\Omega_h$
  as claimed.
\end{proof}
This lemma implies that $\Sigma_h$ minimizes area on the outside.
\begin{lemma}
  \label{thm:outward_minimizing}
  For all $h\in(0,\Hmax)$ and all sets of finite perimeter
  $\Omega_h\subset F\subset \Omega$, we have
  \begin{equation*}
    |\Sigma_h| \leq |\del^+ F|,
  \end{equation*}
  in particular
  \begin{equation*}
    |\Sigma_h|
    \leq
    |\del^+\Omega|.
  \end{equation*}  
\end{lemma}
\begin{proof}
  As $\Omega_h$ minimizes $J_h$ on the outside,
  \begin{equation*}
    |\Sigma_h| + h(\Vol(F)-\Vol(\Omega)) \leq |\del^+ F|.
  \end{equation*}
\end{proof}
To conclude, we mention two other properties, which follow from the
construction.
\begin{lemma}
  \label{thm:basic-properties}
  \begin{enumerate}
  \item\label{item:1} The sets $\Omega_h$ are increasing, that is
    $\Omega_{h_1}\subset \Omega_{h_2}$ if $h_1< h_2$.
  \item\label{item:2}
    If $h\in(0,\Hmax)$ is fixed and $h_k\in (0,\Hmax)$ a sequence
    with $h_k \geq h$ and $\lim_k h_k = h$, then
    \begin{equation*}
      \bar\Omega_h = \bigcap_{k\geq 1} \Omega_{h_k}.
    \end{equation*}
    Here $\bar\Omega_h$ denotes the closure of $\Omega_h$ in $\Omega$.
  \end{enumerate}
\end{lemma}
\begin{proof}
  Property~\ref{item:1} follows from theorem~\ref{thm:exist_barrier},
  as we can always use $\Sigma_{h_1}$ and $\del^+\Omega$ as inner and
  outer barriers for the construction of a surface with CMC $h_2$
  outside. Note that this requires the strong maximum principle to
  conclude that $\Sigma_{h_1}$ is disjoint from $\del^+\Omega$.

  To prove property~\ref{item:2}, note that clearly $\bar\Omega_h
  \subset \bigcap_{k\geq 1} \Omega_{h_k}$, as there is a positive
  distance between $\Sigma_h$ and $\Sigma_{h'}$ if $h<h'$. On the
  other hand, in view of the curvature bound on $\Sigma_h$ and the
  area estimate, lemma~\ref{thm:outward_minimizing}, we can assume
  that the $\Sigma_{h_k}$ converge to a smooth surface $\Sigma'$ with
  CMC $h$. By construction, $\Sigma'$ lies on the outside of
  $\Sigma_h$ and hence must agree with $\Sigma_h$, as $\Sigma_h$ is
  outermost.
\end{proof}
\subsection{Level set formulation}
Clearly, the sets $\Omega_h$ constructed above can be recognized as
the sub-level sets of a function $u$. For $x\in\Omega$, we can define
$u(x)$ as follows
\begin{equation*}
  u(x):= \inf\{ h : x\in\Omega_h \}.
\end{equation*}
We denote the sub-level sets by
\begin{equation*}
  \begin{split}
    E_h
    &
    :=
    \{ x\in \Omega: u(x) < h \},\ \text{and}
    \\
    E_h^+
    &
    :=
    \{ x\in \Omega: u(x) \leq h \}.
  \end{split}
\end{equation*}
We can say the following about these level sets.
\begin{lemma}
  For all $h\in[0,\Hmax]$ we have that $E_h\subset \Omega_h$ and
  $E_h^+ = \bar \Omega_h$.
\end{lemma}
\begin{proof}
  If $x\in E_h$ then $u(x)< h$ which implies $x\in\Omega_h$ by the
  definition of $u$, hence $E_h\subset \Omega_h$.

  Let $x\in E_h^+$, that is $u(x)\leq h$. Then for all $h'>h$ we have
  that $x\in\Omega_{h'}$. As the intersection of all $\Omega_{h'}$
  with $h'>h$ is $\bar\Omega_h$ by property~\ref{item:2} of
  lemma~\ref{thm:basic-properties}, we infer that $E_h^+\subset
  \bar\Omega_h$. To see the other inclusion, note that if
  $x\in\bar\Omega_h$ then $x\in\Omega_{h'}$ for all $h'>h$.
\end{proof}
\begin{lemma}
  If $u$ is as defined above, then $u\in BV(\Omega)\cap C^0(\Omega)$,
  where $BV(\Omega)$ denotes the space of functions with bounded
  variation and $C^0(\Omega)$ denotes the space of bounded continuous
  functions.
\end{lemma}
\begin{proof}
  First note that $u(x)\in [0,\Hmax]$ and thus $u$ is bounded.  

  We show that $u$ is continuous. First, note that since $E_h^+$ is
  closed, we have that $\{u>h\} = \Omega\setminus E_h^+$ is
  open. Furthermore
  \begin{equation*}
    \{u=h\} = \bar\Omega_h \setminus \bigcup_{h'<h} \Omega_{h'}
  \end{equation*}
  hence $\{u<h\} \bigcup_{h'<h} \Omega_{h'}$ and thus $\{u<h\}$ is
  also open. These two properties imply the continuity of $u$.

  Furthermore, for all $k\in\IN$ we can choose values $0=h^k_0 <
  \ldots < h^K_{N(k)}=\Hmax$ such that $|h_i^k - h_{i-1}^k|< 1/k$ for
  $i=1,\ldots, N(k)$. Let
  \begin{equation*}
    u_k := \sum_{i=1}^{N(k)} (h^k_i - h^k_{i-1}) \chi_{E^+_{h_i^k}},
  \end{equation*}
  where $\chi_E$ denotes the characteristic function of a set
  $E$. Note that the $u_k$ converge uniformly to $u$ as $k\to\infty$
  since $u$ is continuous. Furthermore, all $u_k$ have their $BV$-norm
  bounded by $|\Sigma_\Hmax| \Hmax$, and thus contain a subsequence
  that converges weakly to a limit $u_\infty \in BV$. As the $u_k$
  converge uniformly to $u$ we have that $u=u_\infty$ and hence $u$ is
  in $BV$ and has $BV$ norm bounded by $|\Sigma_\Hmax| \Hmax$.
\end{proof}
In \cite{Huisken-Ilmanen:2001}, Huisken and Ilmanen introduced a
notion of weak solutions to the level-set inverse mean curvature
flow. This notion motivates the following definition of a
self-referencing functional on sets $F$ of bounded variation.
\begin{equation*}
  J_u(F) := |\del^+ F| - \int_F u dx.
\end{equation*}
Based on this functional we introduce the notion of weak CMC
foliations.
\begin{definition}
  We say that $u$ is a weak (respectively sub-, super-) solution to
  the CMC foliation problem, if the sets $E_h^+ := \{ x\in M : u(x)
  \leq h\}$ minimize $J_u$ (from the outside, inside respectively).
\end{definition}
With  respect to the above definition, we show the following theorem.
\begin{theorem}
  The function $u$, as defined above, is a weak sub-solution to the
  CMC foliation problem.
\end{theorem}
\begin{proof}
  Let $F\supset \bar\Omega_h$ be any subset of finite perimeter. Fix
  $\eps>0$ and pick $h_i \in(0,\Hmax)$ such that
  \begin{equation*}
    h = h_0 < h_1 < \ldots < h_N ,
  \end{equation*}
  $h_i -h_{i-1} <\eps$ and $h_N$ is such that
  $F\subset\Omega_{h_N}$. For each $h_i$ we know that $\Omega_{h_i}$
  minimizes $J_{h_i}$ from the outside. Hence we can compare with the
  set $F_i := (\Omega_{h_{i+1}} \cap F) \cup \bar\Omega_{h_i}$ and find
  that
  \begin{equation*}
    J_{h_i}(\Omega_{h_i}) \leq J_{h_i}(F_i).
  \end{equation*}
  Expanding this out, we obtain
  \begin{equation*}
    \begin{split}
      |\Sigma_{h_i}| - h_i \Vol(\Omega_{h_i})
      &
      \leq
      |\del^+ F \cap (\bar\Omega_{h_{i+1}} \setminus \bar\Omega_{h_i})|
      + | \Sigma_{h_{i+1}} \cap F |
      + | \Sigma_{h_i} \setminus F |
      \\
      &\phantom{\leq}
      - h_{i} \Vol(\Omega_{h_i}) - h_i \Vol(F\cap(\bar\Omega_{h_{i+1}}
      \setminus \bar \Omega_{h_i}))
    \end{split}
  \end{equation*}
  sorting terms, this implies that
  \begin{equation*}
    |\Sigma_{h_i}\cap F| - |\Sigma_{h_{i+1}}\cap F|
    \leq
    |\del^+ F \cap (\bar\Omega_{h_{i+1}} \setminus \bar\Omega_{h_i})|
    - h_i \Vol(F\cap(\bar\Omega_{h_{i+1}} \setminus \bar \Omega_{h_i})).
  \end{equation*}
  Taking the sum, we find that
  \begin{equation*}
    \sum_{i=0}^{N-1}
    ( |\Sigma_{h_i}\cap F| - |\Sigma_{h_{i+1}}\cap F| )
    \leq
    |\del^+ F|
    - \sum_{i=0}^{N-1}
    h_i\Vol(F\cap(\bar\Omega_{h_{i+1}} \setminus \bar \Omega_{h_i})).
  \end{equation*}
  Since $\Omega_h \subset F$ and $F \subset \bar\Omega_{h_N}$ we have
  that $|\Sigma_{h_0}\cap F| = |\Sigma_h|$ and $|\Sigma_{h_N}\cap F| =
  0$. So the above implies
  \begin{equation}
    \label{eq:6}
    |\Sigma_h| \leq |\del^+ F| - \sum_{i=0}^{N-1}
    h_i\Vol(F\cap(\bar\Omega_{h_{i+1}} \setminus \bar \Omega_{h_i})).
  \end{equation}
  As $h_{i}\leq u\leq h_{i+1}$ on $\Omega_{h_{i+1}}\setminus
  \Omega_{h_i}$ we can estimate
  \begin{equation*}
    \begin{split}
      \int_F u dx - \int_{\Omega_h} u dx
      &=
      \sum_{i=0}^{N-1} \int_{F\cap (\Omega_{h_{i+1}}\setminus \Omega_{h_i})} u dx
      \\
      &\leq
      \sum_{i=0}^{N-1} h_{i+1} \Vol(F\cap (\Omega_{h_{i+1}}\setminus
      \Omega_{h_i}))
      \\
      &\leq
      \sum_{i=0}^{N-1} (h_{i}+\eps) \Vol(F\cap (\Omega_{h_{i+1}}\setminus
      \Omega_{h_i}))
      \\
      &\leq
      \sum_{i=0}^{N-1} h_{i} \Vol(F\cap (\Omega_{h_{i+1}}\setminus
      \Omega_{h_i}))
      + \eps \Vol(F\setminus \Omega_h).      
    \end{split}
  \end{equation*}
  Combining this estimate with equation~\eqref{eq:6} from above, we
  arrive at
  \begin{equation*}
    \int_F u dx - \int_{\Omega_h} u dx
    \leq
    |\del^+ F| - |\Sigma_h| + \eps \Vol(F\setminus \Omega_h)
  \end{equation*}
  This implies that
  \begin{equation*}
    J_u (\Omega_h) \leq J_u (F) + \eps \Vol(F\setminus \Omega_h)
  \end{equation*}
  as $\eps$ was arbitrary, this yields the claim.
\end{proof}
\begin{remark}
  We arrived at a weak sub-solution to the CMC foliation problem in the
  interior region by taking the outermost sets with curvature
  $\Omega_h$. Analogously, we can construct the sets $\tilde \Omega_h$
  bounded by the innermost surfaces with CMC $h$. Then the procedure
  above will result in surfaces minimizing $J_h$ from the inside,
  which in turn implies that the corresponding level set function
  $\tilde u$ is a super-solution to the weak CMC foliation problem.

  Having these sub- and super- solutions at hand it should be possible
  to construct a weak solution of the CMC foliation problem in the
  sense as defined above. This is research in progress, details of
  which will appear elsewhere.
\end{remark}


%
%
%
\section*{Acknowledgements}
Research on this project started while the author was visiting Stanford
University. During this time, the author also received financial
support from the Alexander-von-Humboldt foundation. 

%
\bibliographystyle{amsalpha}
\bibliography{extern/references}

\providecommand{\bysame}{\leavevmode\hbox to3em{\hrulefill}\thinspace}
\providecommand{\MR}{\relax\ifhmode\unskip\space\fi MR }
\providecommand{\MRhref}[2]{%
  \href{http://www.ams.org/mathscinet-getitem?mr=#1}{#2}
}
\providecommand{\href}[2]{#2}
\begin{thebibliography}{SSY75}

\bibitem[AM05]{Andersson-Metzger:2005}
L.~Andersson and J.~Metzger, \emph{Curvature estimates for stable marginally
  trapped surfaces}, arXiv:gr-qc/0512106, 2005.

\bibitem[AM07]{Andersson-Metzger:2007}
\bysame, \emph{The area of horizons and the trapped region}, arXiv:0708.4252
  [gr-qc], 2007.

\bibitem[Bra97]{Bray:1997}
H.~L. Bray, \emph{The {P}enrose inequality in {G}eneral {R}elativity and volume
  comparison theorems involving scalar curvat ure}, Ph.D. thesis, Stanford
  University, 1997.

\bibitem[Gal06]{Galloway:2006ws}
Gregory~J. Galloway, \emph{Rigidity of outer horizons and the topology of black
  holes}, arXiv:gr-qc/0608118, 2006.

\bibitem[HI01]{Huisken-Ilmanen:2001}
G.~Huisken and T.~Ilmanen, \emph{The inverse mean curvature flow and the
  {R}iemannian {P}enrose inequality}, J. Differential Geom. \textbf{59} (2001),
  no.~3, 353--437.

\bibitem[HY96]{Huisken-Yau:1996}
G.~Huisken and S.-T. Yau, \emph{Definition of center of mass for isolated
  physical systems and unique foliations by stable spheres with constant mean
  curvature}, Invent. Math. \textbf{124} (1996), no.~1-3, 281--311.

\bibitem[Pog81]{Pogorelov:1981}
A.~V. Pogorelov, \emph{On the stability of minimal surfaces}, Dokl. Akad. Nauk
  SSSR \textbf{260} (1981), no.~2, 293--295.

\bibitem[Ros06]{Rosenberg:2006}
Harold Rosenberg, \emph{Constant mean curvature surfaces in homogeneously
  regular 3-manifolds}, Bull. Austral. Math. Soc. \textbf{74} (2006), no.~2,
  227--238. \MR{MR2260491 (2007g:53009)}

\bibitem[SSY75]{Schoen-Simon-Yau:1975}
R.~Schoen, L.~Simon, and S.~T. Yau, \emph{Curvature estimates for minimal
  hypersurfaces}, Acta Math. \textbf{134} (1975), no.~3-4, 275--288.

\bibitem[Ye91]{ye:1991}
Rugang Ye, \emph{Foliation by constant mean curvature spheres}, Pacific J.
  Math. \textbf{147} (1991), no.~2, 381--396. \MR{MR1084717 (92f:53030)}

\bibitem[ZZ98]{Ziemer-Zumbrum:1998}
W.~P. Ziemer and K.~Zumbrum, \emph{The obstacle problem for functions of least
  gradient}, Mathematica Bohemica \textbf{124} (1998), 193--219.

\end{thebibliography}
\end{document}